\documentclass[11pt, reqno]{amsart}
\usepackage{amsmath, amsthm, amscd, amsfonts, amssymb, graphicx, color}
\usepackage[bookmarksnumbered, colorlinks, plainpages]{hyperref}
\hypersetup{colorlinks=true,linkcolor=red, anchorcolor=green, citecolor=cyan, urlcolor=red, filecolor=magenta, pdftoolbar=true}

\newtheorem{theorem}{Theorem}[section]
\newtheorem{lemma}[theorem]{Lemma}
\newtheorem*{theorem*}{Question}
\newtheorem{proposition}[theorem]{Proposition}
\newtheorem{corollary}[theorem]{Corollary}
\theoremstyle{definition}
\newtheorem{definition}[theorem]{Definition}
\newtheorem{example}[theorem]{Example}

\theoremstyle{remark}
\newtheorem{remark}[theorem]{Remark}
\numberwithin{equation}{section}

\begin{document}

\setcounter{page}{1}


\title[On the Convexity of Berezin range and Berezin Radius Inequalities] {On the Convexity of Berezin range and Berezin Radius Inequalities via a Class of Semi-Norms}
\author[Athul Augustine, P. Hiran Das, Pintu Bhunia \MakeLowercase{and} P. Shankar]{Athul Augustine, P. Hiran Das, Pintu Bhunia \MakeLowercase{and} P. Shankar}

\address{Athul Augustine, Department of Mathematics, Cochin University of Science And Technology,  Ernakulam, Kerala - 682022, India. }
\email{\textcolor[rgb]{0.00,0.00,0.84}{athulaugus@gmail.com, athulaugus@cusat.ac.in}}

\address{P. Hiran Das, Department of Mathematics, Cochin University of Science And Technology,  Ernakulam, Kerala - 682022, India. }
\email{\textcolor[rgb]{0.00,0.00,0.84}{hirandas939@gmail.com, hirandas073@cusat.ac.in}}

\address{Pintu Bhunia, Department of Mathematics, Indian Institute of Science, Bengaluru 560012, Karnataka, India and Department of Mathematics, SRM University AP, Amaravati 522240, Andhra Pradesh, India.}
\email{\textcolor[rgb]{0.00,0.00,0.84}{pintubhunia5206@gmail.com}}

\address{P. Shankar, Department of Mathematics, Cochin University of Science And Technology,  Ernakulam, Kerala - 682022, India.}
\email{\textcolor[rgb]{0.00,0.00,0.84}{shankarsupy@gmail.com, shankarsupy@cusat.ac.in}}




\begin{abstract}
 This paper introduces a new family of semi-norms, say  $\sigma_\mu$-Berezin norm on the space of all bounded linear operators $B(\mathcal{H})$ defined on a reproducing kernel Hilbert space $\mathcal{H}$, namely, for each $\mu \in [0,1]$ and $p\geq 1$,  
$$\|T\|_{\sigma_{\mu}\text{-ber}}= \sup_{\lambda\in\Omega}\left\lbrace \left(|\langle T\hat{k}_\lambda,\hat{k}_\lambda\rangle |^p~ \sigma_{\mu}~ \|T\hat{k}_\lambda\|^p\right)^{\frac{1}{p}}\right\rbrace $$
where $T\in B(\mathcal{H})$ and $\sigma_{\mu}$ is an interpolation path of the symmetric mean $\sigma$. We investigate many fundamental properties of the $\sigma_\mu$-Berezin norm and develop several inequalities associated with it. Utilizing these inequalities, we derive improved bounds for the Berezin radius of bounded linear operators, enhancing previously known estimates. Furthermore, we study the convexity of the Berezin range of a class of composition operators and weighted shift operators on both the Hardy space and the Bergman space.
\end{abstract}
\maketitle

\section{Introduction}
Let $\Omega$ be a non empty set. A \textit{reproducing kernel Hilbert space} (RKHS) is a Hilbert space of complex-valued functions defined on $\Omega$, such that for every $w \in \Omega$, the evaluation functional $f \mapsto f(w)$ is bounded. By the Riez representation theorem, if $\mathcal{H}$ is a Hilbert space on $\Omega$, for each $\lambda\in\Omega$ and every $f\in \mathcal{H}$, there exist a unique element $k_\lambda$ satisfying, $E_\lambda(f) = f(\lambda)=\langle f,k_\lambda\rangle$. The element $k_\lambda$ is the reproducing kernel at $\lambda$ and we denote the normalized reproducing kernel at $\lambda$ as $\hat{k}_\lambda = \frac{k_\lambda}{\|k_\lambda\|}$. Further details on reproducing kernel Hilbert space can be found in \cite{paulsen2016introduction}.

Let $B(\mathcal{H})$ denote the $C^*$-algebra of all bounded linear operators on the reproducing kernel Hilbert space $\mathcal{H}$. For a bounded linear operator $T\in B(\mathcal{H})$, the Berezin transform (or Berezin symbol) of $T$ at $\lambda$ is defined by $\widetilde{T}(\lambda):= \langle T \hat{k}_\lambda, \hat{k}_\lambda \rangle$. The Berezin range (or Berezin set) and the Berezin radius of an operator $T$ are defined, respectively as
$$ \textit{Ber}(T) := \{\widetilde{T}(\lambda) : \lambda \in \Omega\}\quad \text{and}\quad \textit{ber}(T) := \sup_{\lambda \in \Omega}\{|\widetilde{T}(\lambda)|\}.$$
 The Berezin transform of an operator was first introduced by F. A. Berezin in \cite{berezin1972covariant}. One can easily observe that the Berezin range of an operator in a reproducing kernel Hilbert space lies within its numerical range and  the Berezin radius of an operator $T$ does not exceed its norm, i.e., $ \textit{ber}(T)\leq \|T\|$. Unlike the numerical range, the Berezin range is not necessarily convex.  The convexity of the Berezin range has been investigated in several works, including \cite{augustine2023composition, augus2024, cowen22, karaev2013reproducing}.
The Berezin norm of $T\in B(\mathcal{H})$ is defined by
$$\|T\|_{ber} := \sup_{\lambda\in \Omega}\left\|T\hat{k}_\lambda\right\|.$$
The Berezin norm, $\|\cdot\|_{ber}$, defines an operator norm on $B(\mathcal{H})$. Moreover, it is straightforward that
$$\textit{ber}(T) \leq \|T\|_{\textit{ber}} \leq \|T\|.
$$

Next we discuss mean functions.
A mean is a non-negative function $\sigma:[0,\infty)\times [0,\infty)\rightarrow [0,\infty)$ that satisfies the following conditions (see \cite{bhatia, geometric}):
\begin{enumerate}
\item[(i)]$\sigma(a,b)\geq 0;$
\item[(ii)]if $a\leq b$, $a\leq \sigma(a,b)\leq b;$
\item[(iii)]$\sigma(a,b)$ is monotone increasing in both a and b;
\item[(iv)]For $\alpha>0$, $\sigma(\alpha a,\alpha b)=\alpha\sigma(a,b)$ (homogeneity);
\item[(v)] $\sigma(a,b)$ is continuous.
\end{enumerate}

For convenience, throughout this paper, we denote $\sigma(a, b)$ by $a \sigma b$. A mean, $\sigma$ is called a symmetric mean, if $a\sigma b= b\sigma a$.  For a symmetric mean $\sigma$, a parameterized operator mean $\sigma_\mu$, for each $\mu\in[0,1]$, is called an interpolation path \cite[Section 5.3]{mond} for $\sigma$ if it satisfies the following conditions:
\begin{enumerate}
\item[(i)]$a\sigma_0 b=b,a\sigma_1 b=a,a \sigma_\frac{1}{2} b=a\sigma b~ \forall ~ a,b\geq0;$
\item[(ii)]$(a\sigma_\mu b)\sigma(a\sigma_\nu b) = a\sigma_\frac{\mu+\nu}{2} b, \forall~\mu,\nu \in [0,1];$
\item[(iii)]For each $0\leq \mu \leq 1$, $\sigma_\mu$ is increasing in each of its components.
\end{enumerate}

Well known examples of symmetric means are the arithmetic mean, the geometric mean and the harmonic mean. The arithmetic mean of two non-negative numbers $a,b$ is defined by $a\nabla b=\frac{a+b}{2}$ and the interpolation path for $\nabla$ is defined by
$$ a\nabla_\mu b= \mu a + (1-\mu)b,\quad 0\leq \mu\leq 1.$$
The geometric mean and the harmonic mean of two non-negative numbers $a,b$ are defined respectively as $a\sharp b =\sqrt{ab}$ and $a ! b= \left(\frac{a^{-1}+b^{-1}}{2}\right)^{-1}$. The interpolation paths for these symmetric means are defined by
$$a\sharp_\mu b = a^\mu b^{1-\mu} ~~ \text{and}~~ a !_\mu b= \left(\mu a^{-1} +(1-\mu) b^{-1}\right)^{-1}, ~~ \text{where} ~0\leq \mu\leq 1.$$
The study of means has been extended from non-negative numbers to bounded operators on Hilbert spaces. An interpolation path of the arithmetic mean that lies between the Berezin radius and the Berezin norm was introduced in \cite{newnorm}.
\begin{definition}
Let $T\in B(\mathcal{H})$ and let $\nabla_{\alpha}$ be an interpolation path of the arithmetic mean $\nabla$. Then for each $\alpha\in [0,1]$, $\alpha$-Berezin norm is defined as
$$\|T\|_{\alpha\text{-ber}}= \sup_{\lambda\in\Omega}\left\lbrace \sqrt{|\langle T\hat{k}_\lambda,\hat{k}_\lambda\rangle |^2 \nabla_{\alpha} \|T\hat{k}_\lambda\|^2}\right\rbrace. $$
\end{definition}
It is straightforward to verify that for $\alpha \in [0,1)$, the $\alpha$-Berezin norm defines a norm on \(B(\mathcal{H})\). When \(\alpha = 1\), it is a norm provided the functional Hilbert space satisfies the ``Ber'' property, i.e.,  for \(T, S \in B(\mathcal{H})\),   $\widetilde{T}(\lambda) = \widetilde{S}(\lambda)$ for all \(\lambda \in \Omega\)  implies \(T = S\). Consequently, the \(\alpha\)-Berezin norm defines a norm for all \(\alpha \in [0,1]\) in many well-known functional Hilbert spaces, including the Hardy space and the Bergman spaces.

In this work, we extend this definition to general  interpolation paths of symmetric means. We introduce the $\sigma_\mu$-Berezin norm on $B(\mathcal{H})$, a new numerical quantity that lie between the Berezin radius and the Berezin norm.
\begin{definition}\label{def}
Let $T\in B(\mathcal{H})$ and let $\sigma_{\mu}$ be an interpolation path of the symmetric mean $\sigma$. Then for $p\geq 1$ and $\mu \in [0,1]$, we define the $\sigma_\mu$-Berezin norm on $B(\mathcal{H})$ as
$$\|T\|_{\sigma_{\mu}\text{-ber}}= \sup_{\lambda\in\Omega}\left\lbrace \left(|\langle T\hat{k}_\lambda,\hat{k}_\lambda\rangle |^p~ \sigma_{\mu}~ \|T\hat{k}_\lambda\|^p\right)^{\frac{1}{p}}\right\rbrace. $$
\end{definition}

 In this paper, we examine the fundamental properties of the $\sigma_\mu$-Berezin norm and derive several improved bounds for the Berezin radius of bounded linear operators in $B(\mathcal{H})$ by establishing new inequalities for the $\sigma_{\mu}$-Berezin norm.  Further, we develop upper bounds for the Berezin radius of sums of product operators. We also examine the convexity of the Berezin range associated with composition operators and a class of weighted shift operators  acting on the classical Hardy space and Bergman space. These results unifies and generalizes several previously known results concerning the convexity of the Berezin range \cite{augustine2023composition, cowen22} and Berezin radius inequalities \cite{newfam, newnorm, upperbound}.

\section{The $\sigma_\mu$-Berezin norm and Berezin radius inequalities}

We begin by examining the fundamental properties of the $\sigma_\mu$-Berezin norm and deriving several analytical bounds for it. These bounds are employed to obtain improved upper estimates for the Berezin radius of bounded linear operators on functional Hilbert spaces.

It is clear from the definition that the $\sigma_\mu$-Berezin norm is a semi-norm and satisfies the relation
\begin{equation}\label{ber}
\textit{ber}(T) \leq \|T\|_{\sigma_{\mu}\text{-ber}} \leq \|T\|_{ber},~\forall ~ \mu\in[0,1].
\end{equation}
It is easy to observe that there exist operators $T \in B(\mathcal{H})$ for which $\|T\|_{\sigma_{\mu}\text{-ber}} \neq \|T^*\|_{\sigma_{\mu}\text{-ber}}$ whenever $\mu \neq 1$, where $T^*$ denotes the adjoint of $T$.
\begin{example}
Consider the functional Hilbert space $\mathcal{H}=\mathbb{C}^2$ and $T=\begin{pmatrix}
1& 1\\
2& 0
\end{pmatrix}$.
If $\sigma_\mu$ is the interpolation path of the geometric mean, then $\|T\|_{\sigma_{\mu}\text{-ber}}= \sqrt{5}^\mu \neq \sqrt{2}^\mu = \|T^*\|_{\sigma_{\mu}\text{-ber}}$ for every $\mu\in(0,1]$.
\end{example}
The following example demonstrates a strict inequality between $\textit{ber}(\cdot ),$ $\|\cdot\|_{\sigma_{\mu}\text{-ber}}$ and $\|\cdot\|_{ber}$.
\begin{example}
Consider the rank one operator $A(f)= \langle f,z\rangle z$ on the classical Hardy space $H^2(\mathbb{D})$. Then we have $\textit{ber}(A) = \frac{1}{4}$ and $\|A\|_{ber}=\frac{1}{2}$. If $\sigma_\mu$ is the interpolation path of the geometric mean, then $\|A\|_{\sigma_{\mu}\text{-ber}}= \sup_{\lambda\in\mathbb{D}}\left\lbrace(1-|\lambda|^2)^{(1-\frac{\mu}{2})}|\lambda|^{(2-\mu)}\right\rbrace=\frac{1}{2}^{(2-\mu)}$ and $\textit{ber}(A)< \|A\|_{\sigma_{\mu}\text{-ber}} < \|A\|_{ber}$ for all $\mu\in (0,1)$.
\end{example}
The first proposition follows immediately from the definition of the $\sigma_\mu$-Berezin norm.
\begin{proposition}\label{prop2.1}
Let $T\in B(\mathcal{H})$ and let $\sigma_{\mu}$ be an interpolation path of the symmetric mean $\sigma$.  Then the following results hold:
\begin{enumerate}
\item[$(i)$] $\|T\|_{\sigma_{0}\text{-ber}} = \textit{ber}(T)$, $\|T\|_{\sigma_{1}\text{-ber}} = \|T\|_{ber}$.
\item[$(ii)$] For each $0\leq \mu \leq 1$, $\|T\|_{\sigma_{\mu}\text{-ber}} \leq \left(\textit{ber}(T)^p\sigma_{\mu} \|T\|_{ber}^p\right)^{\frac{1}{p}}$.
\item[$(iii)$] $\|T\|_{\sigma_{\mu}\text{-ber}} = 0$, for any $\mu$, if and only if $T=0$.
\item[$(iv)$] $\|\lambda T\|_{\sigma_{\mu}\text{-ber}} = |\lambda| \|T\|_{\sigma_{\mu}\text{-ber}}$ for any $\lambda \in \mathbb{C}$.
\end{enumerate}
\end{proposition}
For $T\in B(\mathcal{H})$, we have $\||T|\|_{\textit{ber}}=\|T\|_{\textit{ber}}$. For a semi-hyponormal operator $T$, the mixed Cauchy-Schwarz inequality $|\langle T \hat{k}_\lambda, \hat{k}_\lambda\rangle|^2\leq  \langle |T| \hat{k}_\lambda, \hat{k}_\lambda\rangle\langle |T^*| \hat{k}_\lambda, \hat{k}_\lambda\rangle$ implies that $\textit{ber}(T)\leq \textit{ber}(|T|)$. In the following result, we examine possible relations involving the norm $\|\cdot\|_{\sigma_{\mu}\text{-ber}}$ in this context.

\begin{theorem}
Let $T\in B(\mathcal{H})$ and let $\sigma_{\mu}$ be an interpolation path of the symmetric mean $\sigma$.  Then the following results hold:
\begin{enumerate}
\item[$(i)$]If $T$ is hyponormal, i.e., $TT^*\leq T^*T$, then $\|T^*\|_{\sigma_{\mu}\text{-ber}} \leq \|T\|_{\sigma_{\mu}\text{-ber}}$.
\item[$(ii)$] If $T$ is co-hyponormal, i.e., $T^*T\leq TT^*$, then $\|T\|_{\sigma_{\mu}\text{-ber}} \leq \|T^*\|_{\sigma_{\mu}\text{-ber}}$.
\item[$(iii)$] If $T$ is semi-hyponormal, i.e., $|T^*|\leq |T|$, then $\|T\|_{\sigma_{\mu}\text{-ber}} \leq \||T|\|_{\sigma_{\mu}\text{-ber}}$.
\item[$(iv)$] If $T$ is $(\alpha,\beta)$-normal, i.e., $\alpha T^*T\leq TT^*\leq \beta T^*T$ for some positive real numbers $\alpha$ and $\beta$ with $\alpha\leq 1\leq \beta$, then $$\alpha\|T\|_{\sigma_{\mu}\text{-ber}}\leq\|T^*\|_{\sigma_{\mu}\text{-ber}} \leq \beta\|T\|_{\sigma_{\mu}\text{-ber}}.$$
\item[$(v)$] If $T$ is normal, then $\|T^*\|_{\sigma_{\mu}\text{-ber}} = \|T\|_{\sigma_{\mu}\text{-ber}}$.
\end{enumerate}
\end{theorem}
\begin{proof}
$(i)$ By hyponormality of $T$, we have $\|T^*\hat{k}_\lambda\|\leq \|T\hat{k}_\lambda\|$ for all $\lambda\in \Omega.$ Also utilizing the basic fact $|\langle T\hat{k}_\lambda,\hat{k}_\lambda\rangle|=|\langle T^*\hat{k}_\lambda,\hat{k}_\lambda\rangle|$  for all $\lambda\in \Omega$ and by the monotonicity of $\sigma_{\mu}$, we conclude that
$$|\langle T^*\hat{k}_\lambda,\hat{k}_\lambda\rangle |^p~ \sigma_{\mu}~ \|T^*\hat{k}_\lambda\|^p\leq |\langle T\hat{k}_\lambda,\hat{k}_\lambda\rangle |^p~ \sigma_{\mu}~ \|T\hat{k}_\lambda\|^p~~\text{for all}~~  \lambda\in \Omega.$$
 Thus, by taking the supremum over all $\lambda\in \Omega$, we get $\|T^*\|_{\sigma_{\mu}\text{-ber}} \leq \|T\|_{\sigma_{\mu}\text{-ber}}$.\\
$(ii)$ By co-hyponormality of $T$, we have $\|T\hat{k}_\lambda\|\leq \|T^*\hat{k}_\lambda\|$ for all $\lambda\in \Omega.$ Then by similar proof in $(i)$, we get the desired result.\\
$(iii)$ From semi-hyponormality and the mixed Cauchy-Schwarz inequality, we have $|\langle T\hat{k}_\lambda,\hat{k}_\lambda\rangle |^p\leq \langle|T|\hat{k}_\lambda,\hat{k}_\lambda\rangle^p$  for all $\lambda\in \Omega.$  Also utilizing the basic fact $\| T\hat{k}_\lambda\|=\| |T|\hat{k}_\lambda\|$  for all $\lambda\in \Omega$ and by the monotonicity of $\sigma_{\mu}$, we conclude that
$$|\langle T\hat{k}_\lambda,\hat{k}_\lambda\rangle |^p~ \sigma_{\mu}~ \|T\hat{k}_\lambda\|^p\leq \langle |T|\hat{k}_\lambda,\hat{k}_\lambda\rangle^p~ \sigma_{\mu}~ \||T|\hat{k}_\lambda\|^p~~\text{for all}~~  \lambda\in \Omega.$$
 Thus, by taking the supremum over all $\lambda\in \Omega$, we get $\|T\|_{\sigma_{\mu}\text{-ber}} \leq \||T|\|_{\sigma_{\mu}\text{-ber}}$.\\
$(iv)$ Since $T$ is  $(\alpha,\beta)$-normal, we have $\alpha\|T\hat{k}_\lambda\|\leq \|T^*\hat{k}_\lambda\|\leq \beta\|T\hat{k}_\lambda\|$ for all $\lambda\in \Omega$. Then by monotonicity of $\sigma_{\mu}$, for all $\lambda\in \Omega$, we have
$$\alpha^p|\langle T\hat{k}_\lambda,\hat{k}_\lambda\rangle |^p~ \sigma_{\mu}~ \alpha^p\|T\hat{k}_\lambda\|^p\leq \langle T^*\hat{k}_\lambda,\hat{k}_\lambda\rangle^p~ \sigma_{\mu}~ \|T^*\hat{k}_\lambda\|^p\leq \beta^p\langle T\hat{k}_\lambda,\hat{k}_\lambda\rangle^p~ \sigma_{\mu}~ \beta^p\|T\hat{k}_\lambda\|^p.$$
Taking the  supremum over all $\lambda\in \Omega$, we get the desired result.\\
$(v)$ Since $T$ is normal it is both hyponormanl and co-hyponormanl. Therefore, from $(i)$ and $(ii)$, we have $\|T^*\|_{\sigma_{\mu}\text{-ber}} = \|T\|_{\sigma_{\mu}\text{-ber}}$.
\end{proof}

To prove our next results, we need the following basic lemmas.

\begin{lemma}\cite{kato}
If $T\in B(\mathcal{H})$ and $0\leq \alpha\leq 1$, then 
\begin{equation}\label{1.4}
|\langle Tx,y\rangle|^2 \leq \langle|T|^{2\alpha} x,x\rangle\left\langle|T^*|^{2(1-\alpha)}y,y\right\rangle,
\end{equation}
for all $x,y\in\mathcal{H}$. In particular, for $\alpha=\frac{1}{2}$,
\begin{equation}\label{1.5}
|\langle Tx,y\rangle|^2 \leq \langle|T| x,x\rangle\left\langle|T^*|y,y\right\rangle,
\end{equation}
for all $x,y\in\mathcal{H}$.
\end{lemma} 
\begin{lemma}\cite[Theorem 5]{notes}
Let $T,S \in \mathcal{B}(\mathcal{H})$ be such that $|T|S=S^*|T|$, and let $\phi$ and $\psi$ be two non-negative continuous functions defined on $[0,\infty)$ such that $\phi(t)\psi(t)=t$ for every $t\geq 0$. Then
\begin{equation}\label{1.6}
|\langle TS x,y\rangle|^2 \leq r(S)\langle \phi^2(|T|)x,x\rangle \langle \psi^2(|T^*|)y,y\rangle,
\end{equation}
for every $x,y\in \mathcal{H}$, where $r(S)$ denotes the spectral radius of $S$.
In particular, if $S=I,$ the identity operator in $\mathcal{B}(\mathcal{H})$, we obtain 
\begin{equation}\label{1.7}
|\langle Tx,y\rangle|^2
\leq \langle \phi^2(|T|)x,x\rangle \langle \psi^2(|T^*|)y,y\rangle,
\end{equation}
for every $x,y\in \mathcal{H}$.
\end{lemma}
\begin{lemma}\cite{simon}
If $T\in B(\mathcal{H})$ is positive, then 
\begin{equation}\label{1.8}
\langle Tx,x\rangle^p \leq \langle T^p x,x\rangle,~ \forall p\geq 1,~ \forall x\in \mathcal{H},~ \|x\|=1.
\end{equation}
The inequality \eqref{1.8} is reversed when $0\leq p\leq 1$.
\end{lemma}
\begin{lemma}\cite{buzano}
If $x,y,e\in\mathcal{H}$ with $\|e\|=1$, then
\begin{equation}\label{1.9}
|\langle x,e\rangle\langle e,y\rangle|\leq \frac{1}{2}(\|x\|\|y\|+|\langle x,y\rangle|).
\end{equation}
\end{lemma}
\begin{lemma}\cite{complex}\label{sum}
For each $i=1,2,...,n$, let $c_i$ be a positive real number. Then
$$\left(\sum_{i=1}^n c_i\right)^p\leq n^{p-1}\sum_{i=1}^n c^p_i,$$
for all $p\geq 1$.
\end{lemma}
A well-known bound for the Berezin radius of operators in a reproducing kernel Hilbert space is given by the following theorem. In this paper, we present improvements to this bound. 
\begin{theorem}\cite[Corollary 3.4]{upperbound}\label{sigma}
Let $A_i,B_i,X_i\in B(\mathcal{H})$ $(i=1,2,...,n)$ and $0<\alpha<1$. Then for $p\geq 1$,
$$\textit{ber}^p\left(\sum_{i=1}^n A_i^*X_iB_i\right)\leq \frac{n^{p-1}}{2}\textit{ber}\left(\sum_{i=1}^n \left(\left[A_i^*|X_i^*|^{2(1-\alpha)}A_i\right]^p+\left[B_i^*|X_i|^{2\alpha}B_i\right]^p\right)\right).$$
\end{theorem}
The following are some special cases of Theorem \ref{sigma}.
\begin{theorem}\cite[Corollary 3.5]{upperbound}\label{bound}
Let $A,B,X\in B(\mathcal{H})$. Then the following results hold:
\begin{itemize}
\item[(i)]$\textit{ber}^p(A)\leq \frac{1}{2}\textit{ber}\left(|A|^p+|A^*|^p\right)\quad \forall~p\geq 1.$
\item[(ii)]$\textit{ber}(A^*XB)\leq \frac{1}{2}\textit{ber}\left(|A|^2+|B|^2\right).$
\item[(iii)]$\textit{ber}(A^*B)\leq \frac{1}{2}\textit{ber}\left(A^*|X^*|A+B^*|X^*|B\right).$
\end{itemize}
\end{theorem}
For $n=2$, substituting $X_i=I, A_1=B_1=T$ and  $A_2=B_2=T^*$ in Theorem \ref{sigma}, we get the following lemma.
\begin{corollary}\label{p}
Let $T\in B(\mathcal{H})$. Then for $p\geq 1$,
$$\textit{ber}^p\left(|T|^2 + |T^*|^2\right)\leq 2^{p-1}\textit{ber}\left(|T|^{2p} + |T^*|^{2p}\right).$$
\end{corollary}

Now we prove the following theorem.

\begin{theorem}
Let $T\in B(\mathcal{H})$ and let $\mu\in (0,1)$. Then for $p\geq 1$, the following conditions are equivalent:
\begin{itemize}
\item[(1)]$\|T\|^p_{\sigma_{\mu}\text{-ber}}=\textit{ber}^p(T)\sigma_{\mu} \|T\|_{ber}^p$.
\item[(2)] There exists a sequence $\{\lambda_n\} \in \Omega$ such that $\lim\limits_{n \rightarrow \infty} |\langle T\hat{k}_{\lambda_n},\hat{k}_{\lambda_n}\rangle | = \textit{ber}(T)$, and $\lim\limits_{n \rightarrow \infty} \|T\hat{k}_{\lambda_n}\| = \|T\|_{ber}$.
\end{itemize}
\end{theorem}
\begin{proof}
First, we prove that $(1)\implies (2)$. Assume that $(1)$ holds. It is clear that there exists a sequence $\{\lambda_n\}\in\Omega$ such that 
$$\|T\|_{\sigma_{\mu}\text{-ber}}^p=\lim\limits_{n \rightarrow \infty} |\langle T\hat{k}_{\lambda_n},\hat{k}_{\lambda_n}\rangle |^p~\sigma_{\mu}~\lim\limits_{n \rightarrow \infty} \|T\hat{k}_{\lambda_n}\|^p.$$
Since $\{|\langle T\hat{k}_{\lambda_n},\hat{k}_{\lambda_n}\rangle |\}$ and $\{\|T\hat{k}_{\lambda_n}\|\}$ are both bounded sequences of real numbers, there exists a subsequence $\{\lambda_{n_k}\}$ of $\{\lambda_n\}$ such that both the sequences are convergent. Therefore, we have
\begin{equation*}
\begin{split}
\textit{ber}(T)^p\sigma_{\mu} \|T\|_{ber}^p &= \|T\|_{\sigma_{\mu}\text{-ber}}^p\\
&=\lim_{k \rightarrow \infty} |\langle T\hat{k}_{\lambda_{n_k}},\hat{k}_{\lambda_{n_k}}\rangle |^p~\sigma_{\mu}~\lim_{k \longrightarrow \infty} \|T\hat{k}_{\lambda_{n_k}}\|^p\\
&\leq \textit{ber}(T)^p\sigma_{\mu} \|T\|_{ber}^p.
\end{split}
\end{equation*}
This implies that $\lim_{k \longrightarrow \infty} |\langle T\hat{k}_{\lambda_{n_k}},\hat{k}_{\lambda_{n_k}}\rangle | = \textit{ber}(T)$, and $\lim_{k \longrightarrow \infty} \|T\hat{k}_{\lambda_{n_k}}\| = \|T\|_{ber}$.
Now, to prove $(2)\implies (1)$, assume $(2)$ holds. We have
\begin{equation*}
\begin{split}
 \|T\|_{\sigma_{\mu}\text{-ber}}^p &= \sup_{\lambda\in\Omega}\left\lbrace |\langle T\hat{k}_\lambda,\hat{k}_\lambda\rangle |^p~ \sigma_{\mu}~ \|T\hat{k}_\lambda\|^p\right\rbrace\\
&\geq \lim\limits_{n \rightarrow \infty} |\langle T\hat{k}_{\lambda_n},\hat{k}_{\lambda_n}\rangle |^p~\sigma_{\mu}~\lim\limits_{n \rightarrow \infty} \|T\hat{k}_{\lambda_n}\|^p\\
&= \textit{ber}(T)^p~\sigma_{\mu}~ \|T\|_{ber}^p.
\end{split}
\end{equation*}
This completes the proof.
\end{proof} 
For $T\in B(\mathcal{H})$, we define $\widetilde{c}(T)$ by the distance from the origin to the Berezin range of $T$. That is, $\widetilde{c}(T) := \inf_{\lambda\in\Omega}|\langle T\hat{k}_{\lambda},\hat{k}_{\lambda}\rangle |$.
In the following theorem, we obtain a lower bound for the $\sigma_\mu$-Berezin norm in terms of $\widetilde{c}(\cdot)$. 
 \begin{theorem}
Let $T\in B(\mathcal{H})$. Then  for $p\geq 1$,
$$\|T\|_{\sigma_{\mu}\text{-ber}}^p\geq \max\left\lbrace\textit{ber}^p(T)~\sigma_\mu ~ \widetilde{c}^{\frac{p}{2}}(T^*T),~\widetilde{c}^p(T)~\sigma_\mu~ \|T\|_{ber}^p\right\rbrace.$$
 \end{theorem}
 \begin{proof}
 We have 
 \begin{equation*}
\begin{split}
 \|T\|_{\sigma_{\mu}\text{-ber}}^p &\geq  |\langle T\hat{k}_\lambda,\hat{k}_\lambda\rangle |^p~ \sigma_{\mu}~ \|T\hat{k}_\lambda\|^p\\
&=|\langle T\hat{k}_\lambda,\hat{k}_\lambda\rangle |^p~ \sigma_{\mu}~ \langle T^*T\hat{k}_\lambda,\hat{k}_\lambda\rangle^\frac{p}{2}\\
&\geq |\langle T\hat{k}_\lambda,\hat{k}_\lambda\rangle |^p~ \sigma_{\mu}~\widetilde{c}^{\frac{p}{2}}(T^*T).
\end{split}
\end{equation*}
Taking the supremum over $\lambda\in\Omega$, we get
\begin{equation}\label{craw1}
\|T\|_{\sigma_{\mu}\text{-ber}}^p \geq \textit{ber}^p(T)~\sigma_\mu ~ \widetilde{c}^{\frac{p}{2}}(T^*T).
\end{equation}
Also, by similar computation and substitution, we get
$$\|T\|_{\sigma_{\mu}\text{-ber}}^p \geq  |\langle T\hat{k}_\lambda,\hat{k}_\lambda\rangle |^p~ \sigma_{\mu}~ \|T\hat{k}_\lambda\|^p\geq \widetilde{c}^p(T)~ \sigma_{\mu}~ \|T\hat{k}_\lambda\|^p.$$
Taking the supremum over $\lambda\in\Omega$, we get
\begin{equation}\label{craw2}
\|T\|_{\sigma_{\mu}\text{-ber}}^p \geq \widetilde{c}(T)^p~\sigma_\mu~ \|T\|_{ber}^p.
\end{equation}
 Combining \eqref{craw1} and \eqref{craw2}, we get the desired result.
 \end{proof}
 For $\mu=\frac{1}{2}$, we get the following lower bound of the Berezin number.
 
\begin{corollary}
Let $T\in B(\mathcal{H})$. If $\sigma_\frac{1}{2}$ is the geometric mean, then, for $p\geq 1$,
$$ \max\left\lbrace\sqrt{\textit{ber}^p(T) \widetilde{c}^{\frac{p}{2}}(T^*T)},~\sqrt{\widetilde{c}^p(T)\|T\|_{ber}^p}\right\rbrace\leq \|T\|_{\sigma_{\frac12}\text{-ber}}^p.$$
\end{corollary}

We now establish an upper bound for the $\sigma_\mu$-Berezin norm. 

 \begin{theorem}
 Let $T\in B(\mathcal{H})$ and let $\phi$ and $\psi$ be two non-negative continuous functions defined on $[0,\infty)$ such that $\phi(t)\psi(t)=t$ for every $t\geq 0$. If $\sigma_\mu\leq \nabla_\mu$, then, for $p\geq 1$,
 \begin{itemize}
 \item[(i)]$\|T\|_{\sigma_{\mu}\text{-ber}}^p \leq \textit{ber}\left(\frac{\mu}{2^p}(\phi^2(|T|) + \psi^2(|T^*|))^p+(1-\mu)  |T|^p\right),$
 \item[(ii)]$\textit{ber}^p(T)\leq \inf_{\mu\in [0,1]}\textit{ber}\left(\frac{\mu}{2^p}(\phi^2(|T|) + \psi^2(|T^*|))^p+(1-\mu)  |T|^p\right).$
 \end{itemize}
 \end{theorem}
\begin{proof}
We have 
\begin{equation*}
\begin{split}
 |\langle T\hat{k}_\lambda,\hat{k}_\lambda\rangle|^2 &\leq \langle \phi^2(|T|)\hat{k}_\lambda,\hat{k}_\lambda\rangle \langle \psi^2(|T^*|)\hat{k}_\lambda,\hat{k}_\lambda\rangle\\
 &=\left(\langle \phi^2(|T|)\hat{k}_\lambda,\hat{k}_\lambda\rangle^{\frac{1}{2}} \langle \psi^2(|T^*|)\hat{k}_\lambda,\hat{k}_\lambda\rangle^{\frac{1}{2}}\right)^2\\
 &\leq \frac{1}{4}\left(\langle \phi^2(|T|)\hat{k}_\lambda,\hat{k}_\lambda\rangle+ \langle \psi^2(|T^*|)\hat{k}_\lambda,\hat{k}_\lambda\rangle\right)^2\\
 &=\frac{1}{4}\left(\langle (\phi^2(|T|) + \psi^2(|T^*|))\hat{k}_\lambda,\hat{k}_\lambda\rangle\right)^2\\
 &\leq\frac{1}{4}\left(\langle (\phi^2(|T|) + \psi^2(|T^*|))^2\hat{k}_\lambda,\hat{k}_\lambda\rangle\right),\\
\end{split}
\end{equation*}
where we have used Equation \eqref{1.7} to obtain the first inequality,  the arithmetic-geometric mean inequality to obtain the second inequality and Equation \eqref{1.8} to obtain the third inequality.
Therefore, from \eqref{1.8}, we have 
\begin{equation*}
\begin{split}
 \|T\|_{\sigma_{\mu}\text{-ber}}^p & =  \sup_{\lambda\in\Omega}\left\lbrace |\langle T\hat{k}_\lambda,\hat{k}_\lambda\rangle |^p~ \sigma_{\mu}~ \|T\hat{k}_\lambda\|^p\right\rbrace\\
 &\leq \sup_{\lambda\in\Omega}\left\lbrace \frac{1}{4^\frac{p}{2}}\left\langle (\phi^2(|T|) + \psi^2(|T^*|))^2\hat{k}_\lambda,\hat{k}_\lambda\right\rangle^\frac{p}{2}~ \sigma_{\mu}~ \langle |T|^2\hat{k}_\lambda,\hat{k}_\lambda\rangle^\frac{p}{2}\right\rbrace\\
 &\leq \sup_{\lambda\in\Omega}\left\lbrace \frac{1}{2^p}\left\langle (\phi^2(|T|) + \psi^2(|T^*|))^p\hat{k}_\lambda,\hat{k}_\lambda\right\rangle~ \sigma_{\mu}~ \langle |T|^p\hat{k}_\lambda,\hat{k}_\lambda\rangle\right\rbrace\\
 &\leq \sup_{\lambda\in\Omega}\left\lbrace \frac{\mu}{2^p}\left\langle (\phi^2(|T|) + \psi^2(|T^*|))^p\hat{k}_\lambda,\hat{k}_\lambda\right\rangle+(1-\mu) \langle |T|^p\hat{k}_\lambda,\hat{k}_\lambda\rangle\right\rbrace\\
 &= \sup_{\lambda\in\Omega}\left\lbrace \left\langle \frac{\mu}{2^p}(\phi^2(|T|) + \psi^2(|T^*|))^p+(1-\mu)  |T|^p\hat{k}_\lambda,\hat{k}_\lambda\right\rangle\right\rbrace\\
  &= \textit{ber}\left(\frac{\mu}{2^p}(\phi^2(|T|) + \psi^2(|T^*|))^p+(1-\mu)  |T|^p\right),
\end{split}
\end{equation*}
as desired.
\end{proof} 
If we take $ \phi(t)= t^\frac{1}{2}$, $\psi(t) =t^\frac{1}{2}$, $\mu=1$ and combining it with \eqref{ber}, we get the following corollary.
\begin{corollary}\label{3.2}
Let $T\in B(\mathcal{H})$. Then, for $p\geq 2$,
$$\textit{ber}^p(T)\leq \frac{1}{2}\textit{ber} \left( |T|^p + |T^*|^p \right).$$
\end{corollary}
\begin{proof}
$$\textit{ber}^p(T) \leq \textit{ber}\left(\frac{1}{2^p}\left(|T| + |T^*|\right)^p\right)\leq \frac{1}{2}\textit{ber} \left( |T|^p + |T^*|^p \right),$$
where we have used Lemma \ref{sum} to obtain the second inequality.
\end{proof}

If we take $ \phi(t)= t^{1-\nu}, \psi(t) =t^\nu$, where $\nu\in(0,1)$, we get the following corollary.
\begin{corollary}\label{1.11}
Let $T\in B(\mathcal{H})$. Then, for $p\geq 1$,
$$\|T\|_{\sigma_{\mu}\text{-ber}}^p \leq 
\textit{ber} \left( \frac{\mu}{2^p} \left( |T|^{2(1 - \nu)} + |T^*|^{2\nu} \right)^p + (1 - \mu) |T|^p \right)$$
and 
$$\textit{ber}^p(T) \leq \inf_{\mu \in [0,1]} \textit{ber} \left( \frac{\mu}{2^p} \left( |T|^{2(1 - \nu)} + |T^*|^{2\nu} \right)^p + (1 - \mu) |T|^p \right).$$
\end{corollary}
The second inequality follows by combining the first part of Corollary \ref{1.11} with Equation \eqref{ber}.

Next result is as follows.
\begin{theorem}
Let $T\in B(\mathcal{H})$. If $\sigma_\mu\leq \nabla_\mu$, then, for $p\geq 1$,
$$\|T\|_{\sigma_{\mu}\text{-ber}}^p \leq \textit{ber}\left(\mu|T^*|^p +(1-\mu)|T|^p\right).$$
\end{theorem}
\begin{proof}
We have
\begin{equation*}
\begin{split}
 \|T\|_{\sigma_{\mu}\text{-ber}}^p & =  \sup_{\lambda\in\Omega}\left\lbrace |\langle T\hat{k}_\lambda,\hat{k}_\lambda\rangle |^p~ \sigma_{\mu}~ \|T\hat{k}_\lambda\|^p\right\rbrace\\
&= \sup_{\lambda\in\Omega}\left\lbrace |\langle \hat{k}_\lambda,T^*\hat{k}_\lambda\rangle |^p~ \sigma_{\mu}~ \|T\hat{k}_\lambda\|^p\right\rbrace\\
 &\leq \sup_{\lambda\in\Omega}\left\lbrace \|T^*\hat{k}_\lambda\|^p~ \sigma_{\mu}~ \|T\hat{k}_\lambda\|^p\right\rbrace\\
&= \sup_{\lambda\in\Omega}\left\lbrace \left\langle|T^*|^2\hat{k}_\lambda,\hat{k}_\lambda\right\rangle^\frac{p}{2}~ \sigma_{\mu}~\left\langle|T|^2\hat{k}_\lambda,\hat{k}_\lambda\right\rangle^\frac{p}{2}\right\rbrace\\ 
&\leq \sup_{\lambda\in\Omega}\left\lbrace \left\langle|T^*|^p\hat{k}_\lambda,\hat{k}_\lambda\right\rangle~ \sigma_{\mu}~\left\langle|T|^p\hat{k}_\lambda,\hat{k}_\lambda\right\rangle\right\rbrace\\
&\leq \sup_{\lambda\in\Omega}\left\lbrace \mu\left\langle|T^*|^p\hat{k}_\lambda,\hat{k}_\lambda\right\rangle +(1-\mu)\left\langle|T|^p\hat{k}_\lambda,\hat{k}_\lambda\right\rangle\right\rbrace\\
&= \sup_{\lambda\in\Omega}\left\lbrace \left\langle \mu |T^*|^p +(1-\mu)|T|^p\hat{k}_\lambda,\hat{k}_\lambda\right\rangle\right\rbrace\\  
  &=\textit{ber}\left(\mu|T^*|^p +(1-\mu)|T|^p\right),    
\end{split}
\end{equation*}
as desired.
\end{proof}
\begin{remark}
For $\mu=\frac{1}{2}$, combining Equation \eqref{ber} with the above theorem, we get $\textit{ber}^p(T)\leq \frac{1}{2}\textit{ber} \left( |T|^p + |T^*|^p \right),$ which was given in \cite[Corollary 3.5 (i)]{upperbound}.
\end{remark}
\begin{corollary}\label{cor3.4}
Let $T\in B(\mathcal{H})$. Then, for $p\geq 1$,
$$\textit{ber}^p(T)\leq \inf_{\mu\in[0,1]}\textit{ber}\left(\mu|T^*|^p +(1-\mu)|T|^p\right)\leq  \frac{1}{2}\textit{ber} \left( |T|^p + |T^*|^p \right).$$
\end{corollary}
\begin{example}
To show the second inequality of Corollary \ref{cor3.4} is proper, we consider the finite-dimensional Hilbert space $\mathbb{C}^3$ and the operator
\[
T = 
\begin{pmatrix}
0 & 3 & 1 \\
0 & 0 & 2 \\
0 & 0 & 0
\end{pmatrix}.
\]
Then, for $p=4$, $\inf_{\mu\in[0,1]}\textit{ber}\left(\mu|T^*|^4 +(1-\mu)|T|^4\right)=\frac{1560}{29}$ and $\frac{1}{2}\textit{ber} \left( |T|^4 + |T^*|^4 \right)=55$. Therefore, for the matrix $T$, we have
$$\inf_{\mu\in[0,1]}\textit{ber}\left(\mu|T^*|^4 +(1-\mu)|T|^4\right)<  \frac{1}{2}\textit{ber} \left( |T|^4 + |T^*|^4 \right).$$
\end{example}
\begin{theorem}\label{thm3.5}
Let $T\in B(\mathcal{H})$.  If $\sigma_\mu\leq \nabla_\mu$, then, for $p\geq 1$,
$$\|T\|_{\sigma_{\mu}\text{-ber}}^p \leq \left( \frac{\mu}{4} \right)^\frac{p}{2} \textit{ber}^{\frac{p}{2}}\left( |T|^2 + |T^*|^2 \right)+2\left( \frac{\mu}{8} \right)^{\frac{p}{2}}\textit{ber}^\frac{p}{2}\left(|T|^2 + |T^*|^2\right) +(1-\mu)\textit{ber}^\frac{p}{2}\left(|T|^2\right).$$
\end{theorem}
\begin{proof}
Applying the inequality \eqref{1.5} and the arithmetic-geometric mean inequality, we get
\begin{equation*}
\begin{split}
 \left|\left\langle T\hat{k}_\lambda,\hat{k}_\lambda\right\rangle\right|^2&\leq \left\langle|T| \hat{k}_\lambda,\hat{k}_\lambda\right\rangle\left\langle|T^*|\hat{k}_\lambda,\hat{k}_\lambda\right\rangle\\
 &\leq \frac{1}{4}\left(\left\langle|T| \hat{k}_\lambda,\hat{k}_\lambda\right\rangle + \left\langle|T^*|\hat{k}_\lambda,\hat{k}_\lambda\right\rangle\right)^2\\
 & = \frac{1}{4}\left(\left\langle|T| \hat{k}_\lambda,\hat{k}_\lambda\right\rangle^2 + \left\langle|T^*|\hat{k}_\lambda,\hat{k}_\lambda\right\rangle^2\right)+\frac{1}{2}\left\langle|T| \hat{k}_\lambda,\hat{k}_\lambda\right\rangle\left\langle|T^*|\hat{k}_\lambda,\hat{k}_\lambda\right\rangle\\
 &=\frac{1}{4}\left|\left\langle|T| \hat{k}_\lambda,\hat{k}_\lambda\right\rangle + i\left\langle|T^*|\hat{k}_\lambda,\hat{k}_\lambda\right\rangle\right|^2+\frac{1}{2}\left\langle|T^*| \hat{k}_\lambda,\hat{k}_\lambda\right\rangle\left\langle\hat{k}_\lambda,|T|\hat{k}_\lambda\right\rangle\\
 &\leq \frac{1}{4}\left|\left\langle (|T| + i|T^*|)\hat{k}_\lambda,\hat{k}_\lambda\right\rangle\right|^2+\frac{1}{4}\||T^*| \hat{k}_\lambda\|\||T| \hat{k}_\lambda\|+\frac{1}{4}\left|\left\langle |T^*|\hat{k}_\lambda,|T|\hat{k}_\lambda\right\rangle\right|\\
 &\leq \frac{1}{4}\left|\left\langle (|T| + i|T^*|)\hat{k}_\lambda,\hat{k}_\lambda\right\rangle\right|^2+\frac{1}{8}\left(\||T^*| \hat{k}_\lambda\|^2+\||T| \hat{k}_\lambda\|^2\right)+\frac{1}{4}\left|\left\langle |T||T^*|\hat{k}_\lambda,\hat{k}_\lambda\right\rangle\right|.\\
 \end{split}
\end{equation*}
Therefore, we obtain
\begin{equation*}
\begin{split}
 \|T\|_{\sigma_{\mu}\text{-ber}}^p & =  \sup_{\lambda\in\Omega}\left\lbrace |\langle T\hat{k}_\lambda,\hat{k}_\lambda\rangle |^p~ \sigma_{\mu}~ \|T\hat{k}_\lambda\|^p\right\rbrace\\
 &\leq \sup_{\lambda\in\Omega} \left(\frac{1}{4}\left|\left\langle (|T| + i|T^*|)\hat{k}_\lambda,\hat{k}_\lambda\right\rangle\right|^2+\frac{1}{8}\left(\||T^*| \hat{k}_\lambda\|^2+\||T| \hat{k}_\lambda\|^2\right)+\frac{1}{4}\left|\left\langle |T||T^*|\hat{k}_\lambda,\hat{k}_\lambda\right\rangle\right|\right)^{\frac{p}{2}}\\
 &\qquad \sigma_{\mu}~\left\langle|T|^2\hat{k}_\lambda,\hat{k}_\lambda\right\rangle^\frac{p}{2}\\
 &\leq \sup_{\lambda\in\Omega} \mu\left(\frac{1}{4}\left|\left\langle (|T| + i|T^*|)\hat{k}_\lambda,\hat{k}_\lambda\right\rangle\right|^2+\frac{1}{8}\left(\||T^*| \hat{k}_\lambda\|^2+\||T| \hat{k}_\lambda\|^2\right)+\frac{1}{4}\left|\left\langle |T||T^*|\hat{k}_\lambda,\hat{k}_\lambda\right\rangle\right|\right)^{\frac{p}{2}}\\
 &\qquad +~(1-\mu)\left\langle|T|^2\hat{k}_\lambda,\hat{k}_\lambda\right\rangle^\frac{p}{2}\\
 &\leq \textit{ber}^p \left(\frac{\sqrt{\mu}}{2}(|T| + i|T^*|\right)+\textit{ber}^\frac{p}{2}\left(\frac{\mu}{8}(|T|^2 + |T^*|^2)\right)+\textit{ber}^\frac{p}{2}\left(\frac{\mu}{4}( |T||T^*|)\right)\\
 &\qquad +(1-\mu)\textit{ber}^\frac{p}{2}\left(|T|^2\right).\\
\end{split}
\end{equation*}
Since $\textit{ber}^2\left(|T| + i|T^*|\right)\leq \textit{ber}\left(|T|^2 + |T^*|^2\right)$ and $\textit{ber}\left(|T||T^*|\right)\leq \frac{1}{2}\textit{ber}\left(|T|^2 + |T^*|^2\right)$, the above inequality become
\begin{equation*}
\begin{split}
 \|T\|_{\sigma_{\mu}\text{-ber}}^p & \leq \left( \frac{\mu}{4} \right)^\frac{p}{2} \textit{ber}^{\frac{p}{2}}\left( |T|^2 + |T^*|^2 \right)+\left( \frac{\mu}{8} \right)^{\frac{p}{2}}\textit{ber}^\frac{p}{2}\left(|T|^2 + |T^*|^2\right)+\left( \frac{\mu}{8} \right)^{\frac{p}{2}} \textit{ber}^{\frac{p}{2}}\left( |T|^2 + |T^*|^2 \right)\\
 &\quad +(1-\mu)\textit{ber}^\frac{p}{2}\left(|T|^2\right)\\
 & = \left( \frac{\mu}{4} \right)^\frac{p}{2} \textit{ber}^{\frac{p}{2}}\left( |T|^2 + |T^*|^2 \right)+2\left( \frac{\mu}{8} \right)^{\frac{p}{2}}\textit{ber}^\frac{p}{2}\left(|T|^2 + |T^*|^2\right) +(1-\mu)\textit{ber}^\frac{p}{2}\left(|T|^2\right),
 \end{split}
\end{equation*}
as desired.
\end{proof}

We deduce the following corollary.

\begin{corollary}\label{cor3.4}
Let $T\in B(\mathcal{H})$.  If $\sigma_\mu\leq \nabla_\mu$, then for $p\geq 2$ and $\mu=1$
$$\textit{ber}^p(T)\leq \left(2^{-p}+2^{-\frac{p}{2}-1}\right) \textit{ber}\left( |T|^p + |T^*|^p \right). $$
\end{corollary}

\begin{proof}
Applying Equation \eqref{ber} and Corollary \ref{p} in Theorem \ref{thm3.5}, we get 
\begin{equation*}
\begin{split}
\textit{ber}^p(T) &\leq \left( \frac{1}{4} \right)^\frac{p}{2} \textit{ber}^{\frac{p}{2}}\left( |T|^2 + |T^*|^2 \right)+2\left( \frac{1}{8} \right)^{\frac{p}{2}}\textit{ber}^\frac{p}{2}\left(|T|^2 + |T^*|^2\right)\\
&\leq \left(2^{-p}+2^{1-\frac{3p}{2}}\right) \textit{ber}^{\frac{p}{2}}\left( |T|^2 + |T^*|^2 \right)\\
&\leq \left(2^{-p}+2^{1-\frac{3p}{2}}\right)2^{\frac{p}{2}-1} \textit{ber}\left( |T|^p + |T^*|^p \right)\\
&\leq \left(2^{-p}+2^{-\frac{p}{2}-1}\right) \textit{ber}\left( |T|^p + |T^*|^p \right).
\end{split}
\end{equation*}
\end{proof}

\begin{remark}
We observe that for $p\geq 2$, Corollary \ref{cor3.4} gives better bound than the existing bound in \cite[Corollary 3.5 (i)]{upperbound}.
\end{remark}

\begin{theorem}
Let $T\in B(\mathcal{H})$. Then, for $p\geq 1$,
$$\|T\|_{\sigma_{\mu}\text{-ber}}^p \leq \sup_{\lambda\in\Omega}\left\lbrace\left\langle \frac{1}{2}\left(|T|^2+ |T^*|^2\right)\hat{k}_\lambda,\hat{k}_\lambda\right\rangle^\frac{p}{2}~ \sigma_{\mu}~\left(\langle |T|^2\hat{k}_\lambda,\hat{k}_\lambda\rangle\right)^\frac{p}{2}\right\rbrace.$$
\end{theorem}

\begin{proof}
We have
\begin{equation*}
\begin{split}
 \|T\|_{\sigma_{\mu}\text{-ber}}^p & =  \sup_{\lambda\in\Omega}\left\lbrace |\langle T\hat{k}_\lambda,\hat{k}_\lambda\rangle |^p~ \sigma_{\mu}~ \|T\hat{k}_\lambda\|^p\right\rbrace\\
 &= \sup_{\lambda\in\Omega}\left\lbrace\left( \langle \mathcal{R}(T)\hat{k}_\lambda,\hat{k}_\lambda\rangle^2+\langle \mathcal{I}(T)\hat{k}_\lambda,\hat{k}_\lambda\rangle^2\right)^\frac{p}{2}~ \sigma_{\mu}~\left(\langle T^*T\hat{k}_\lambda,\hat{k}_\lambda\rangle\right)^\frac{p}{2}\right\rbrace\\
 &\leq \sup_{\lambda\in\Omega}\left\lbrace\left( \| \mathcal{R}(T)\hat{k}_\lambda\|^2+\|\mathcal{I}(T)\hat{k}_\lambda\|^2\right)^\frac{p}{2}~ \sigma_{\mu}~\left(\langle T^*T\hat{k}_\lambda,\hat{k}_\lambda\rangle\right)^\frac{p}{2}\right\rbrace\\
 &= \sup_{\lambda\in\Omega}\left\lbrace\left( \langle \mathcal{R}^2(T)\hat{k}_\lambda,\hat{k}_\lambda\rangle+\langle \mathcal{I}^2(T)\hat{k}_\lambda,\hat{k}_\lambda\rangle\right)^\frac{p}{2}~ \sigma_{\mu}~\left(\langle T^*T\hat{k}_\lambda,\hat{k}_\lambda\rangle\right)^\frac{p}{2}\right\rbrace\\
 &= \sup_{\lambda\in\Omega}\left\lbrace\left( \langle \mathcal{R}^2(T)+ \mathcal{I}^2(T)\hat{k}_\lambda,\hat{k}_\lambda\rangle\right)^\frac{p}{2}~ \sigma_{\mu}~\left(\langle |T|^2\hat{k}_\lambda,\hat{k}_\lambda\rangle\right)^\frac{p}{2}\right\rbrace\\
 &= \sup_{\lambda\in\Omega}\left\lbrace\left\langle \frac{1}{2}\left(|T|^2+ |T^*|^2\right)\hat{k}_\lambda,\hat{k}_\lambda\right\rangle^\frac{p}{2}~ \sigma_{\mu}~\left(\langle |T|^2\hat{k}_\lambda,\hat{k}_\lambda\rangle\right)^\frac{p}{2}\right\rbrace.
\end{split}
\end{equation*}
\end{proof}

\begin{remark}
For $\mu =0$ and $p\geq 2$, by Corollary \ref{p}, we have 
$$\textit{ber}^p(T) \leq \frac{1}{2^{\frac{p}{2}}}\textit{ber}^{\frac{p}{2}}\left(|T|^2+|T^*|^2\right)
\leq \frac{1}{2}\textit{ber} \left( |T|^p + |T^*|^p \right).$$
\end{remark}

We now present the following theorem, which provides a new upper bound for the $\sigma_\mu$-Berezin norm and the Berezin radius.

\begin{theorem}
Let $T\in B(\mathcal{H})$. Then, for $p\geq 1$,
$$\|T\|_{\sigma_{\mu}\text{-ber}}^p \leq \sup_{\lambda\in\Omega}\left\lbrace \left[\frac{1}{4}\left(\left\langle (T^*T+ TT^*)\hat{k}_\lambda,\hat{k}_\lambda\right\rangle\right) +\frac{1}{2} \langle T^2\hat{k}_\lambda,\hat{k}_\lambda\rangle \right]^\frac{p}{2}~ \sigma_{\mu}~ \left(\langle |T|^2\hat{k}_\lambda,\hat{k}_\lambda\rangle\right)^\frac{p}{2}\right\rbrace.$$
\end{theorem}
\begin{proof}
Applying the arithmetic-geometric mean inequality, we get
\begin{equation*}
\begin{split}
|\langle T\hat{k}_\lambda,\hat{k}_\lambda\rangle |^2&\leq \left[\frac{1}{2}\left(\|T\hat{k}_\lambda\|\|T\hat{k}_\lambda\| + |\langle T\hat{k}_\lambda,\hat{k}_\lambda\rangle |^2\right)\right]^\frac{p}{2}\\
&=\left[\frac{1}{2}\langle T^*T\hat{k}_\lambda,\hat{k}_\lambda\rangle^\frac{1}{2}\langle TT^*\hat{k}_\lambda,\hat{k}_\lambda\rangle^\frac{1}{2} +\frac{1}{2} |\langle T\hat{k}_\lambda,\hat{k}_\lambda\rangle |^2\right]^\frac{p}{2}\\
&\leq \left[\frac{1}{4}\left(\langle T^*T\hat{k}_\lambda,\hat{k}_\lambda\rangle+\langle TT^*\hat{k}_\lambda,\hat{k}_\lambda\rangle\right) +\frac{1}{2} |\langle T\hat{k}_\lambda,\hat{k}_\lambda\rangle |^2\right]^\frac{p}{2}\\
&= \left[\frac{1}{4}\left(\left\langle (T^*T+ TT^*)\hat{k}_\lambda,\hat{k}_\lambda\right\rangle\right) +\frac{1}{2} |\langle T\hat{k}_\lambda,\hat{k}_\lambda\rangle |^2\right]^\frac{p}{2}.
\end{split}
\end{equation*}
Hence, we have
\begin{equation*}
\begin{split}
 \|T\|_{\sigma_{\mu}\text{-ber}}^p & =  \sup_{\lambda\in\Omega}\left\lbrace |\langle T\hat{k}_\lambda,\hat{k}_\lambda\rangle |^p~ \sigma_{\mu}~ \|T\hat{k}_\lambda\|^p\right\rbrace\\
 &\leq \sup_{\lambda\in\Omega}\left\lbrace \left[\frac{1}{4}\left(\left\langle (T^*T+ TT^*)\hat{k}_\lambda,\hat{k}_\lambda\right\rangle\right) +\frac{1}{2} \langle T^2\hat{k}_\lambda,\hat{k}_\lambda\rangle \right]^\frac{p}{2}~ \sigma_{\mu}~ \left(\langle |T|^2\hat{k}_\lambda,\hat{k}_\lambda\rangle\right)^\frac{p}{2}\right\rbrace.\\
\end{split}
\end{equation*}
\end{proof}

For $\mu=0$, we get the following upper bound of the Berezin number.

\begin{corollary}\label{3.7}
Let $T\in B(\mathcal{H})$. If $T^2=0$, then for $p\geq 2$
$$\textit{ber}^p(T) \leq \frac{1}{2^{1+\frac{p}{2}}}\textit{ber}(|T|^p + |T^*|^p).$$
\end{corollary}
\begin{proof}
Applying Lemma \ref{sigma}, we have
$$\textit{ber}^p(T) \leq \frac{1}{4^{\frac{p}{2}}}\textit{ber}^{\frac{p}{2}}\left(|T|^2+|T^*|^2\right)
\leq \frac{1}{2^{1+\frac{p}{2}}}\textit{ber} \left( |T|^p + |T^*|^p \right).$$
\end{proof}
\begin{remark}
We observe that for $p\geq 2$, Corollary \ref{3.7} gives better bound than the existing bound in \cite[Corollary 3.5 (i)]{upperbound}.
\end{remark}
\begin{theorem}
Let $T\in B(\mathcal{H})$ and $p\geq 1$. Then
$$\|T\|_{\sigma_{\mu}\text{-ber}}^p \leq \sup_{\lambda\in\Omega}\left\lbrace \left(\langle |\mathcal{R}(T)|+ |\mathcal{I}(T)|\hat{k}_\lambda,\hat{k}_\lambda\rangle \right)^p~ \sigma_{\mu}~ \langle T^*T\hat{k}_\lambda,\hat{k}_\lambda\rangle^\frac{p}{2}\right\rbrace.$$
\end{theorem}

\begin{proof}
We have
\begin{equation*}
\begin{split}
 \|T\|_{\sigma_{\mu}\text{-ber}}^p & =  \sup_{\lambda\in\Omega}\left\lbrace |\langle T\hat{k}_\lambda,\hat{k}_\lambda\rangle |^p~ \sigma_{\mu}~ \|T\hat{k}_\lambda\|^p\right\rbrace\\
 &\leq \sup_{\lambda\in\Omega}\left\lbrace \left|\langle \mathcal{R}(T)\hat{k}_\lambda,\hat{k}_\lambda\rangle+i\langle \mathcal{I}(T)\hat{k}_\lambda,\hat{k}_\lambda\rangle \right|^p~ \sigma_{\mu}~ \langle T^*T\hat{k}_\lambda,\hat{k}_\lambda\rangle^\frac{p}{2}\right\rbrace\\
 &\leq \sup_{\lambda\in\Omega}\left\lbrace \left(|\langle \mathcal{R}(T)\hat{k}_\lambda,\hat{k}_\lambda\rangle |+|\langle \mathcal{I}(T)\hat{k}_\lambda,\hat{k}_\lambda\rangle| \right)^p~ \sigma_{\mu}~ \langle T^*T\hat{k}_\lambda,\hat{k}_\lambda\rangle^\frac{p}{2}\right\rbrace\\
 &\leq \sup_{\lambda\in\Omega}\left\lbrace \left(\langle |\mathcal{R}(T)|\hat{k}_\lambda,\hat{k}_\lambda\rangle +\langle |\mathcal{I}(T)|\hat{k}_\lambda,\hat{k}_\lambda\rangle \right)^p~ \sigma_{\mu}~ \langle T^*T\hat{k}_\lambda,\hat{k}_\lambda\rangle^\frac{p}{2}\right\rbrace\\
 &= \sup_{\lambda\in\Omega}\left\lbrace \left(\langle |\mathcal{R}(T)|+ |\mathcal{I}(T)|\hat{k}_\lambda,\hat{k}_\lambda\rangle \right)^p~ \sigma_{\mu}~ \langle T^*T\hat{k}_\lambda,\hat{k}_\lambda\rangle^\frac{p}{2}\right\rbrace.\\
\end{split}
\end{equation*}
\end{proof}
\begin{theorem}
Let $T,S \in \mathcal{B}(\mathcal{H})$ be such that $|T|S=S^*|T|$, and let $\phi$ and $\psi$ be two non-negative continuous functions defined on $[0,\infty)$ such that $\phi(t)\psi(t)=t$ for every $t\geq 0$. Then, for $p\geq 1$,
$$\|TS\|_{\sigma_{\mu}\text{-ber}}^p \leq \sup_{\lambda\in\Omega}\left\lbrace \frac{r^{\frac{p}{2}}(S)}{2}\left(\langle (\phi^{2p}(|T|)+ \psi^{2p}(|T^*|))\hat{k}_\lambda,\hat{k}_\lambda\rangle\right)~ \sigma_{\mu}~ \langle |TS|^2\hat{k}_\lambda,\hat{k}_\lambda\rangle^\frac{p}{2}\right\rbrace,$$
 where $r(S)$ denotes the spectral radius of $S$.

\end{theorem}
\begin{proof}
Applying Equation \eqref{1.6}, the arithmetic-geometric mean inequality and Lemma \ref{sum}, we get
\begin{equation*}
\begin{split}
 |\langle TS\hat{k}_\lambda,\hat{k}_\lambda\rangle |^p&~ \sigma_{\mu}~ \|TS\hat{k}_\lambda\|^p\\
 &=|\langle TS\hat{k}_\lambda,\hat{k}_\lambda\rangle |^p~ \sigma_{\mu}~ \langle |TS|^2\hat{k}_\lambda,\hat{k}_\lambda\rangle^\frac{p}{2}\\
 &\leq \left(\sqrt{r(S)}\langle \phi^2(|T|)\hat{k}_\lambda,\hat{k}_\lambda\rangle^\frac{1}{2} \langle \psi^2(|T^*|)\hat{k}_\lambda,\hat{k}_\lambda\rangle^\frac{1}{2}\right)^p~ \sigma_{\mu}~ \langle |TS|^2\hat{k}_\lambda,\hat{k}_\lambda\rangle^\frac{p}{2}\\
 &\leq \left(\frac{\sqrt{r(S)}}{2}\left(\langle \phi^2(|T|)\hat{k}_\lambda,\hat{k}_\lambda\rangle +\langle \psi^2(|T^*|)\hat{k}_\lambda,\hat{k}_\lambda\rangle\right)\right)^p~ \sigma_{\mu}~ \langle |TS|^2\hat{k}_\lambda,\hat{k}_\lambda\rangle^\frac{p}{2}\\
  &= \frac{r^{\frac{p}{2}}(S)}{2^{p}}\left(\langle \phi^2(|T|)\hat{k}_\lambda,\hat{k}_\lambda\rangle +\langle \psi^2(|T^*|)\hat{k}_\lambda,\hat{k}_\lambda\rangle\right)^p~ \sigma_{\mu}~ \langle |TS|^2\hat{k}_\lambda,\hat{k}_\lambda\rangle^\frac{p}{2}\\
   &\leq \frac{r^{\frac{p}{2}}(S)}{2}\left(\langle \phi^{2p}(|T|)\hat{k}_\lambda,\hat{k}_\lambda\rangle +\langle \psi^{2p}(|T^*|)\hat{k}_\lambda,\hat{k}_\lambda\rangle\right)~ \sigma_{\mu}~ \langle |TS|^2\hat{k}_\lambda,\hat{k}_\lambda\rangle^\frac{p}{2}\\
&= \frac{r^{\frac{p}{2}}(S)}{2}\left(\langle (\phi^{2p}(|T|)+ \psi^{2p}(|T^*|))\hat{k}_\lambda,\hat{k}_\lambda\rangle\right)~ \sigma_{\mu}~ \langle |TS|^2\hat{k}_\lambda,\hat{k}_\lambda\rangle^\frac{p}{2}.
\end{split}
\end{equation*}
Now, by taking supremum over all $\lambda\in \Omega$, we obtain the desired result.
\end{proof}

Next theorem is

\begin{theorem}\label{2.25}
Let $T_i,S_i \in \mathcal{B}(\mathcal{H})$ be such that $|T_i|S_i=S_i^*|T_i|$ for $i=1,2,...,n$ and let $\phi$ and $\psi$ be two non-negative continuous functions defined on $[0,\infty)$ such that $\phi(t)\psi(t)=t$ for every $t\geq 0$. Then, for all $p\geq 1,$ 
$$\textit{ber}^p\left(\sum_{i=1}^n T_iS_i\right)\leq \frac{n^{p-1}}{2}\textit{ber}\left(\sum_{i=1}^nr^\frac{p}{2}(S_i)\left(\phi^{2p}(|T_i|)+  \psi^{2p}(|T_i^*|)\right)\right),$$
where $r(S_i)$ denotes the spectral radius of $S_i$.
\end{theorem}
\begin{proof}
$\left|\left\langle \left(\sum_{i=1}^n T_iS_i\right)\hat{k}_\lambda,\hat{k}_\lambda\right\rangle \right|^p$
\begin{equation*}
\begin{split}
&= \left|\sum_{i=1}^n\left\langle  T_iS_i\hat{k}_\lambda,\hat{k}_\lambda\right\rangle \right|^p\\
&\leq \left(\sum_{i=1}^n\left|\left\langle  T_iS_i\hat{k}_\lambda,\hat{k}_\lambda\right\rangle \right|\right)^p\\
&\leq \left(\sum_{i=1}^n\left(r(S_i)\langle \phi^2(|T_i|)\hat{k}_\lambda,\hat{k}_\lambda\rangle \langle \psi^2(|T_i^*|)\hat{k}_\lambda,\hat{k}_\lambda\rangle\right)^\frac{1}{2}\right)^p\\
&\leq \left(\sum_{i=1}^n \sqrt{r(S_i)}\frac{\langle \phi^2(|T_i|)\hat{k}_\lambda,\hat{k}_\lambda\rangle+ \langle \psi^2(|T_i^*|)\hat{k}_\lambda,\hat{k}_\lambda\rangle}{2}\right)^p\\
&= \frac{n^{p-1}}{2^p}\sum_{i=1}^n r^\frac{p}{2}(S_i)\left(\langle \phi^2(|T_i|)\hat{k}_\lambda,\hat{k}_\lambda\rangle+ \langle \psi^2(|T_i^*|)\hat{k}_\lambda,\hat{k}_\lambda\rangle\right)^p\\
&\leq \frac{n^{p-1}}{2}\sum_{i=1}^nr^\frac{p}{2}(S_i) \left(\langle \phi^{2p}(|T_i|)\hat{k}_\lambda,\hat{k}_\lambda\rangle+ \langle \psi^{2p}(|T_i^*|)\hat{k}_\lambda,\hat{k}_\lambda\rangle\right)\\
&\leq \frac{n^{p-1}}{2}\textit{ber}\left(\sum_{i=1}^nr^\frac{p}{2}(S_i)\left(\phi^{2p}(|T_i|)+  \psi^{2p}(|T_i^*|)\right)\right).
\end{split}
\end{equation*}
 Now, taking the supremum over all $\lambda\in \Omega$, we obtain the required result.
\end{proof}
Considering $n=1$, $T_1=T$, $S_1=S$ and $\phi(t)=\psi(t)=\sqrt{t}$ in Theorem \ref{2.25}, we obtain the following bounds.
\begin{corollary}
Let $T,S \in \mathcal{B}(\mathcal{H})$ be such that $|T|S=S^*|T|$. Then
\begin{itemize}
\item[(i)]$\textit{ber}^p\left( TS\right)\leq \frac{r^\frac{p}{2}(S)}{2}\textit{ber}\left(|T|^p+ |T^*|^p\right).$
\item[(ii)]$\textit{ber}^p\left( T\right)\leq \frac{1}{2}\textit{ber}\left(|T|^p+ |T^*|^p\right).$
\end{itemize}
\end{corollary}
In particular, if $S_i=I$ for each $i=1,2,...,n$,  we obtain the following inequality for the sum of operators.
\begin{corollary}
Let $T_i \in \mathcal{B}(\mathcal{H})$ be such that $|T_i|S_i=S_i^*|T_i|$ for $i=1,2,...,n$ and let $\phi$ and $\psi$ be two non-negative continuous functions defined on $[0,\infty)$ such that $\phi(t)\psi(t)=t$ for every $t\geq 0$. Then, for all $p\geq 1$ 
$$\textit{ber}^p\left(\sum_{i=1}^n T_i\right)\leq \frac{n^{p-1}}{2}\textit{ber}\left(\sum_{i=1}^n\left(\phi^{2p}(|T_i|)+  \psi^{2p}(|T_i^*|)\right)\right).$$
\end{corollary}
By setting $n = 1$ and $f(t) = g(t) = \sqrt{t}$ in the above Corollary, we get the following upper bound for the Berezin radius, which was also given in \cite[Corollary 3.5 (i)]{upperbound}.
\begin{corollary}
Let $T \in \mathcal{B}(\mathcal{H})$. Then, for all $p\geq 1$ 
$$\textit{ber}^p\left(T\right)\leq \frac{1}{2}\textit{ber}\left(|T|^p+  |T^*|^p\right).$$
\end{corollary}

\section{Convexity of the Berezin range}

This section is devoted to analyzing the convexity of the Berezin range for certain composition operators and weighted shifts defined on the Hardy and Bergman spaces. Studies on the convexity of the Berezin range of composition operators can be found in \cite{augustine2023composition, cowen22}. This section presents a generalization of the convexity result for composition operators with elliptic symbols, specifically considering the function $\phi(z) = \zeta |z|^k z$ with $k \geq 0$. We further introduce weighted shift operators on the Hardy and Bergman spaces and study the convexity of the associated Berezin range under the weight sequence $\beta_n = \beta^{n-1}$ for some fixed $\beta \in \mathbb{D}$.

\subsection{On Hardy space}
Let $\mathbb{D}$ be the open unit disc and Hol$(\mathbb{D})$ denotes the set of all holomorphic functions in $\mathbb{D}$. The classical \textit{Hilbert Hardy space} on the unit disc $\mathbb{D}$, is defined as 
$$H^{2}(\mathbb{D}) := \left\lbrace f(z) = \displaystyle\sum_{n\geq 0}a_{n}z^{n} \in \text{Hol}(\mathbb{D}) :\displaystyle\sum_{n\geq 0}|a_{n}|^{2} < \infty \right\rbrace.$$
 Let $f(z)=\sum_{n\geq 0}a_{n}z^{n}$ and $g(z)=\sum_{n\geq 0}b_{n}z^{n}$ be the elements of the Hardy space $H^2(\mathbb{D})$. The inner product between $f$ and $g$ is defined as $\langle f,g \rangle = \sum_{n\geq 0}a_{n}\overline{b_n}$. The space $H^{2}(\mathbb{D})$ is a reproducing kernel Hilbert space and its reproducing kernel,  known as the Szego kernel, is given by
 $$k_{w}(z) = \frac{1}{1 - \bar{w}z},\quad z,w \in \mathbb{D}.$$
 Furthermore, $\|k_\lambda\|^2=\langle k_\lambda, k_\lambda\rangle=k_\lambda(\lambda)=\frac{1}{1 - |\lambda|^2}$.  Let $\phi$ be a complex-valued function $\phi : \mathbb{D}\longrightarrow \mathbb{D}$. A composition operator $C_\phi$ acting on $H^{2}(\mathbb{D})$ is defined by $ C_\phi f := f \circ \phi.$ Then the Berezin transform of $C_\phi$ at $\lambda$ is
$$\widetilde{C_{\phi}}(\lambda) =\langle C_{\phi}\hat{k}_{\lambda},\hat{k}_{\lambda}\rangle=(1-|\lambda|^2)\langle C_{\phi}{k_{\lambda}},{k_{\lambda}}\rangle=(1-|\lambda|^2)k_{\lambda}(\phi(\lambda)).$$
In general, the Berezin range of a composition operator on $H^{2}(\mathbb{D})$ is not convex. Characterizations of the convexity of the Berezin range of $C_\phi$ on the Hardy space can be found in \cite{augustine2023composition, cowen22}. Here we prove a more extended result.

\begin{theorem}
Let $\zeta \in \overline{\mathbb{D}}$ and $\phi(z)=\zeta|z|^kz$ where $k\geq 0$ be an analytic function on the unit disc $\mathbb{D}$. Then the Berezin range of the composition operator $C_{\phi}$ acting on $H^2(\mathbb{D})$ is convex if and only if $\zeta \in [-1,1]$.
\end{theorem}
\begin{proof}
Suppose that the Berezin range of $C_{\phi}$ is convex. Put $z=re^{i\theta}$. Then, we have 
$$\widetilde{C}_\phi (re^{i\theta})= \frac{1-r^2}{1- r^{k+2}\zeta}.$$
Since the function is independent of $\theta$, $\textit{Ber}(C_\phi)$ is a path in $\mathbb{C}$. The assumption of convexity will imply that $\textit{Ber}(C_\phi)$ is either a point or a line segment. $\textit{Ber}(C_\phi)$ is a point if and only if $\zeta =1$ and $k=0$. So, assume that $\textit{Ber}(C_\phi)$  is a line segment. For $z=0$, we have $\widetilde{C}_\phi (0)= 1$ and when taking the radial limit, $\lim_{r\rightarrow 1^-}\widetilde{C}_\phi (re^{i\theta})= 0$. Therefore $\textit{Ber}(C_\phi)$ must be a line 
segment passing through the point $1$ and approaching $0$. This implies that $\textit{Ber}(C_\phi)$ is real and therefore the imaginary part of $\textit{Ber}(C_\phi)$ is $0$. This happens if and only if imaginary part of $\zeta$ is zero. Therefore $\zeta \in [-1,1]$.

Conversely, assume that $\zeta \in [-1,1]$. If $\zeta =1$ and $k=0$, then $\textit{Ber}(C_\phi) = \{1\}$, which is convex in $\mathbb{C}$. Now for $k\geq 0$, we have $\widetilde{C}_\phi (re^{i\theta})= \frac{1-r^2}{1- r^{k+2}\zeta}$. Since $r^{k+2}\zeta \leq r^2$ and $\lim_{r\rightarrow 1^-}\widetilde{C}_\phi (re^{i\theta})= 0$, we have
$$\textit{Ber}(C_\phi) = \left\lbrace \frac{1-r^2}{1- r^{k+2}\zeta} : r\in[0,1)\right\rbrace = (0,1],$$
which is also a convex set.
\end{proof}

\begin{remark}
When $k=0$, \cite[Theorem 4.1]{augustine2023composition} and \cite[Theorem 4.1]{cowen22} can be viewed as corollaries of the above theorem.
\end{remark}

Now, we characterise the convexity of the Berezin range for a class of weighted shift operator on $H^2(\mathbb{D})$. Let $f(z)=\sum_{n\geq 0} a_n z^n$ is an element of $H^2(\mathbb{D})$. We define the weighted shift operator on $H^2(\mathbb{D})$ as
$$T\left(\sum_{n\geq 0} a_n z^n\right) := \sum_{n\geq 0} a_n \beta_{n+1}z^{n+1},$$
where the weights $\{\beta_n\}_{n\geq 0}$ is a bounded sequence in $\mathbb{C}$. Since $\{\beta_n\}_{n\geq 0}$ is bounded, there exist a real number $M$ such that $|\beta_n|<M$ for all $n\geq 0$. Therefore
$$\sum_{n\geq 0} |a_n \beta_{n+1}|^2 < \sum_{n\geq 0} M^2|a_n|^2 < \infty.$$
Hence $T(f)\in H^2(\mathbb{D})$ and $T$ is well defined. The Berezin transform
\begin{equation*}
\begin{split}
\widetilde{T}(\lambda) &= \langle T\hat{k}_\lambda, \hat{k}_\lambda\rangle\\
&= (1-|\lambda|^2)\left\langle \sum_{n\geq 0} \overline{\lambda}^n \beta_{n+1}z^{n+1},  \sum_{n\geq 0} \overline{\lambda}^n z^n\right\rangle\\
&= (1-|\lambda|^2)\lambda \left(\sum_{n\geq 0} |\lambda|^{2n}\beta_{n+1}\right).
\end{split}
\end{equation*}
Since $|\lambda|<1$ and $|\beta_n|<M$, we have 
\begin{equation*}
\begin{split}
|\widetilde{T}(\lambda)| &= \left|(1-|\lambda|^2)\lambda \left(\sum_{n\geq 0} |\lambda|^{2n}\beta_{n+1}\right)\right|\\
&\leq (1-|\lambda|^2)|\lambda| M \left(\sum_{n\geq 0} |\lambda|^{2n}\right)\\
&< M.
\end{split}
\end{equation*}
Therefore, we have
\begin{equation}\label{range}
\textit{Ber}(T) \subseteq \{z\in\mathbb{C}: |z|<M\}
\end{equation}
and 
\begin{equation}\label{radi}
\textit{ber}(T) <M.
\end{equation}
If $\beta_n=c$ for all $n$, then $\widetilde{T}(\lambda)=c\lambda$ and therefore, $\textit{Ber}(T) = \{z\in\mathbb{C}: |z|<|c|\}$ and $\textit{ber}(T) = |c|.$
\begin{theorem}
Let $T$ be the weighted shift operator on $H^2(\mathbb{D})$, with weights given by $\beta_n = \beta^{n-1}$ for some fixed $\beta\in \mathbb{D}$. Then the Berezin range of $T$ is a disc centred at the origin and is therefore convex.
\end{theorem}
\begin{proof}
Substituting the weights on the Berezin transform, we have
\begin{equation*}
\begin{split}
\widetilde{T}(\lambda) &=(1-|\lambda|^2)\lambda \left(\sum_{n\geq 0} |\lambda|^{2n}\beta^n\right)\\
&= \frac{(1-|\lambda|^2)\lambda}{1-\beta|\lambda|^2}.
\end{split}
\end{equation*}
Put $\lambda=re^{i\theta}$ and let $\beta=|\beta|e^{i\psi}$. Then we have 
\begin{equation*}
\begin{split}
\widetilde{T}(re^{i\theta}) &= \frac{(1-r^2)re^{i\theta}}{1-r^2|\beta|e^{i\psi}}\\
&=\frac{(1-r^2)r}{1-r^2|\beta|e^{i\psi}}e^{i\theta}.
\end{split}
\end{equation*}
Therefore, $\textit{Ber}(T)$ is disc centred at the origin and radius
\begin{equation}\label{radius}
\textit{ber}(T)=\sup_{r\in[0,1)} \left|\frac{(1-r^2)r}{1-r^2|\beta|e^{i\psi}}\right| \leq\sup_{r\in[0,1)}\frac{(1-r^2)r}{1-r^2|\beta|}.
\end{equation}
\end{proof}

\begin{remark}
In order to find the Berezin radius of $T$, we try to find the critical points of Equation \eqref{radius}. But this expression becomes quite complicated to solve in general, unless we make some assumptions on $\beta$ to simplify the analysis. So we try to find an upper bound for the Berezin radius of $T$. The numerator $(1-r^2)r$ attains its maximum at $r=\sqrt{\frac{1}{3}}$ and the denominator attains its minimum when $r\rightarrow 1^-$. Therefore,
$$\textit{ber}(T)\leq\sup_{r\in[0,1)} \frac{(1-r^2)r}{1-r^2|\beta|} \leq \frac{2}{3\sqrt{3}(1-|\beta|)}.$$
\end{remark}

\subsection{On Bergman space}
Let $\mathbb{D}$ be the open unit disc and Hol$(\mathbb{D})$ denotes the set of all holomorphic functions in $\mathbb{D}$. The \textit{Bergman space} \cite{paulsen2016introduction} on the unit disc $\mathbb{D}$ is defined by
 $$A^2(\mathbb{D}) = \left\lbrace f\in \text{Hol}(\mathbb{D}): \int_{\mathbb{D}} |f(z)|^2dV(z) < \infty\right\rbrace, $$
 where $dV$ is the normalized area measure on $\mathbb{D}$. The space $A^{2}(\mathbb{D})$ is a reproducing kernel Hilbert space and  the reproducing kernel for $A^2(\mathbb{D})$ is 
 $$ k_w(z) = \frac{1}{(1-\bar{w}z)^2},\quad z,w \in \mathbb{D}.$$
Also, note that  $\|k_\lambda\|^2=\langle k_\lambda, k_\lambda\rangle=k_\lambda(\lambda)=\frac{1}{(1 - |\lambda|^2)^2}$.  Let $\phi$ be a complex-valued function $\phi : \mathbb{D}\longrightarrow \mathbb{D}$. A composition operator $C_\phi$ acting on $A^{2}(\mathbb{D})$ is defined by $ C_\phi f := f \circ \phi.$ Then the Berezin transform of $C_\phi$ at $\lambda$ is
$$\widetilde{C_{\phi}}(\lambda) =\langle C_{\phi}\hat{k}_{\lambda},\hat{k}_{\lambda}\rangle
  			=(1-|\lambda|^2)^2\langle C_{\phi}{k_{\lambda}},{k_{\lambda}}\rangle=(1-|\lambda|^2)^2k_{\lambda}(\phi(\lambda)).$$
In general, the Berezin range of a composition operator on $A^{2}(\mathbb{D})$ is not convex. Characterization of the convexity of the Berezin range of $C_\phi$ on the Bergman space can be found in \cite{augustine2023composition}. Here we obtain a more extended result.

\begin{theorem}\label{00.1}
Let $\zeta \in \overline{\mathbb{D}}$ and $\phi(z)=\zeta|z|^kz$ where $k\geq 0$ be an analytic function on the unit disc $\mathbb{D}$. Then the Berezin range of the composition operator $C_{\phi}$ acting on $H^2(\mathbb{D})$ is convex if and only if $\zeta \in [-1,1]$.
\end{theorem}
\begin{proof}
Suppose that the Berezin range of $C_{\phi}$ is convex. Put $z=re^{i\theta}$. Then, we have 
$$\widetilde{C}_\phi (re^{i\theta})= \left[\frac{(1-r^2)}{(1- r^{k+2}\zeta)}\right]^2.$$
Since the function is independent of $\theta$, $\textit{Ber}(C_\phi)$ is a path in $\mathbb{C}$. The assumption of convexity will imply that $\textit{Ber}(C_\phi)$ is either a point or a line segment. $\textit{Ber}(C_\phi)$ is a point if and only if $\zeta =1$ and $k=0$. So, assume that $\textit{Ber}(C_\phi)$  is a line segment. For $z=0$, we have $\widetilde{C}_\phi (0)= 1$ and when taking the radial limit, $\lim_{r\rightarrow 1^-}\widetilde{C}_\phi (re^{i\theta})= 0$. Therefore $\textit{Ber}(C_\phi)$ must be a line 
segment passing through the point $1$ and approaching $0$. This implies that $\textit{Ber}(C_\phi)$ is real and therefore the imaginary part of $\textit{Ber}(C_\phi)$ is $0$. This happens if and only if imaginary part of $\zeta$ is zero. Therefore $\zeta \in [-1,1]$.

Conversely, assume that $\zeta \in [-1,1]$. If $\zeta =1$ and $k=0$, then $\textit{Ber}(C_\phi) = \{1\}$, which is convex in $\mathbb{C}$. Now for $k\geq 0$, we have $\widetilde{C}_\phi (re^{i\theta})= \frac{1-r^2}{1- r^{k+2}\zeta}$. Since $r^{k+2}\zeta \leq r^2$ and $\lim_{r\rightarrow 1^-}\widetilde{C}_\phi (re^{i\theta})= 0$, we have
$$\textit{Ber}(C_\phi) = \left\lbrace \left[\frac{(1-r^2)}{(1- r^{k+2}\zeta)}\right]^2 : r\in[0,1)\right\rbrace = (0,1],$$
which is also a convex set.
\end{proof}

\begin{remark}
For the case $k=0$, \cite[Theorem 5.1]{augustine2023composition}
can be viewed as a corollary of Theorem \ref{00.1}.
\end{remark}

Now, we characterise the convexity of the Berezin range for a class of weighted shift operators on the Bergman space.
$A^2(\mathbb{D})$ has the orthonormal basis $\{e_n\}_{n=0}^{\infty}$, where 
$$e_n(z) = \sqrt{n+1}z^n \text{ for all } n\geq 0.$$
For $f,g\in A^2(\mathbb{D})$  the inner product can also be expressed as
$$\langle f,g\rangle = \sum_{n=0}^{\infty} \frac{1}{n+1}a_n\overline{b_n},$$
where $f(z)=\sum_{n=0}^{\infty} a_nz^n$ and $g(z)=\sum_{n=0}^{\infty} b_nz^n$.
We have $$k_\lambda(z) = \sum_{n=0}^{\infty}(n+1) \overline{\lambda}^nz^n.$$
Let $f(z)=\sum_{n\geq 0} a_n z^n$ is an element of $A^2(\mathbb{D})$. We define the weighted shift operator on $A^2(\mathbb{D})$ as
$$T\left(\sum_{n\geq 0} a_n \sqrt{n+1}z^n\right) := \sum_{n\geq 0} a_n \beta_{n+1}z^{n+1},$$
where the weights $\{\beta_n\}_{n\geq 0}$ is a bounded sequence in $\mathbb{C}$. Since $\{\beta_n\}_{n\geq 0}$ is bounded, there exist a real number $M$ such that $|\beta_n|<M$ for all $n\geq 0$. Therefore
$$
\|Tf\|^2=\langle Tf,Tf\rangle = \sum_{n=0}^{\infty} \frac{|a_n|^2|b_n|^2}{n+1}
\leq M^2\sum_{n=0}^{\infty}\frac{|a_n|^2}{n+1}< \infty.
$$
Hence $T(f)\in A^2(\mathbb{D})$ and $T$ is well defined. The Berezin transform

\begin{equation*}
\begin{split}
\widetilde{T}(\lambda) &= \langle T\hat{k}_\lambda, \hat{k}_\lambda\rangle\\
&= (1-|\lambda|^2)^2\left\langle \sum_{n\geq 0} (n+1)\overline{\lambda}^n \beta_{n+1}z^{n+1},  \sum_{n\geq 0} (n+1)\overline{\lambda}^n z^n\right\rangle\\
&= (1-|\lambda|^2)^2\lambda \left(\sum_{n\geq 0} (n+1)|\lambda|^{2n}\beta_{n+1}\right).
\end{split}
\end{equation*}
Since $|\lambda|<1$ and $|\beta_n|<M$, we have 
\begin{equation*}
\begin{split}
|\widetilde{T}(\lambda)| &= \left|(1-|\lambda|^2)^2\lambda \left(\sum_{n\geq 0} (n+1)|\lambda|^{2n}\beta_{n+1}\right)\right|\\
&\leq (1-|\lambda|^2)^2|\lambda| M \left(\sum_{n\geq 0} (n+1)|\lambda|^{2n}\right)\\
&< M.
\end{split}
\end{equation*}
Therefore, we have
\begin{equation}\label{bergrange}
\textit{Ber}(T) \subseteq \{z\in\mathbb{C}: |z|<M\}
\end{equation}
and 
\begin{equation}\label{bergradi}
\textit{ber}(T) <M.
\end{equation}
If $\beta_n=c$ for all $n$, then $\widetilde{T}(\lambda)=c\lambda$ and therefore, $\textit{Ber}(T) = \{z\in\mathbb{C}: |z|<|c|\}$ and $\textit{ber}(T) = |c|.$
\begin{theorem}
Let $T$ be the weighted shift operator on $A^2(\mathbb{D})$, with weights given by $\beta_n = \beta^{n-1}$ for some fixed $\beta\in \mathbb{D}$. Then the Berezin range of $T$ is a disc centred at the origin and is therefore convex.
\end{theorem}
\begin{proof}
Substituting the weights on the Berezin transform, we have
\begin{equation*}
\begin{split}
\widetilde{T}(\lambda) &=(1-|\lambda|^2)^2\lambda \left(\sum_{n\geq 0} (n+1)|\lambda|^{2n}\beta^n\right)\\
&= \frac{(1-|\lambda|^2)^2\lambda}{1-\beta|\lambda|^2}.
\end{split}
\end{equation*}
Put $\lambda=re^{i\theta}$ and let $\beta=|\beta|e^{i\psi}$. Then we have 
\begin{equation*}
\begin{split}
\widetilde{T}(re^{i\theta}) &= \frac{(1-r^2)^2re^{i\theta}}{1-r^2|\beta|e^{i\psi}}\\
&=\frac{(1-r^2)^2r}{1-r^2|\beta|e^{i\psi}}e^{i\theta}.
\end{split}
\end{equation*}
Therefore, $\textit{Ber}(T)$ is disc centred at the origin and radius
\begin{equation}\label{bergmanradius}
\textit{ber}(T)=\sup_{r\in[0,1)} \left|\frac{(1-r^2)^2r}{1-r^2|\beta|e^{i\psi}}\right| \leq\sup_{r\in[0,1)}\frac{(1-r^2)^2r}{1-r^2|\beta|}.
\end{equation}
\end{proof}
\begin{remark}
In order to find the Berezin radius of $T$, we try to find the critical points of Equation \eqref{bergmanradius}. But this expression becomes quite complicated to solve in general, unless we make some assumptions on $\beta$ to simplify the analysis. So we try to find an upper bound for the Berezin radius of $T$. The numerator $(1-r^2)^2r$ attains its maximum at $r=\sqrt{\frac{1}{5}}$ and the denominator attains its minimum when $r\rightarrow 1^-$. Therefore
$$\textit{ber}(T)\leq\sup_{r\in[0,1)} \frac{(1-r^2)^2r}{1-r^2|\beta|} \leq \frac{16}{25\sqrt{5}(1-|\beta|)}.$$
\end{remark}

\textbf{Declaration of competing interest}

There is no competing interest.\\

\textbf{Data availability}

No data was used for the research described in the article.\\

{\bf Acknowledgments.} The first author is supported by the Senior Research Fellowship (09/0239(13298)/2022-EMR-I) of CSIR (Council of Scientific and Industrial Research, India). The second author is supported by the Junior Research Fellowship of UGC (University Grants Commission, India). The third author was supported by National Post-Doctoral Fellowship PDF/2022/000325 from SERB (Govt.\ of India), SwarnaJayanti Fellowship SB/SJF/2019-20/14  (PI: Apoorva Khare) from SERB (Govt. of India) and NBHM Post-Doctoral Fellowship 0204/16(3)/2024/R\&D-II/6747 from National Board for Higher Mathematics (Govt.\ of India). The fourth author is supported by the Teachers Association for Research Excellence (TARE/2022/000063) of SERB (Science and Engineering Research Board, India).

\nocite{*}
\bibliographystyle{amsplain}

\end{document}